\documentclass[leqno,a4paper]{article}

\input{tcilatex}
\begin{document}

\title{Flow with Nonlinear Potential in General Networks -- Simulation,
Optimization, Control, Risk and Stability Analysis}
\author{Emmanuel M. Livshits, Leonid A. Ostromuhov\thanks{%
Dr. L. A. Ostromuhov, Wingas, Friedrich-Ebert-Str 160, D-34119 Kassel,
Germany.}\thanks{\textit{e-mail:} \textit{leonid.ostromuhov@wingas.de}}%
\thanks{%
Presented on 16th IMACS World Congress 2000 on Scientific Computation,
Applied Mathematics and Simulation in Lausanne, Switzerland. -- IMACS World
Congresses on Computational and Applied Mathematics and Applications in
Science and Engineering}}
\date{}
\maketitle

\begin{abstract}
The aim of this paper is a short survey of models and methods that developed
by the authors. These models and methods are used to optimize general
networks with nonlinear non-convex restrictions and objectives possessing
mixed continuous-discrete optimization variables. There are discussed the
problem formulations and solution methods for simulation, optimization,
sensitivity and stability analysis for flow with nonlinear potential in
general networks. These problems and the developed methods and programs have
industrial application e.g. by gas networks.
\end{abstract}

\tableofcontents


\markright{Flow with Nonlinear Potential in General Networks}





\section{Introduction}

There are 2 approaches to describe flow in real and virtual networks
presently. The first approach is based on a physical description of flow
laws. Conservation laws of mass, impulse and energy resulting in Ohm law for
electrical current or Bernoulli law for fluid flow, and thermic state
equation are used there. This approach led to physical and technical
modelling and simulation of real existing systems. Dynamic optimization
problem is not solved. Instead of dynamic optimization, 'try and error'
method based on simulation was proposed.

The second approach is developed more for economics. It is based on graph
theory and special methods of linear programming. The flow conservation law
represents a main part of restrictions in a mathematical model there. Such
mathematical problems as maximal flow in network and minimal cost network
flow could be mentioned there. In their dual problems, linear potential is
appeared. Linear potential is connected with flow by means of linear
analogue of Ohm law for passive electric chains.

However, there exist technical systems within which the connection between
flow and potential is not linear. The objective function there depends not
only either on flow or on potential, but on both of them. There are
non-convex and non-smooth functions among constraints and objectives. The
pipe networks, energy systems and electrical chains with nonlinear
associated facilities are examples for such large - scale technical systems.
There are some working approaches to optimize such systems \cite{LO-84}, 
\cite{Mur-86}, \cite{O-98} but the optimization theory based on dynamic
models is not developed for general networks with nonlinear connection
between flow and potential.

The aim of this paper is a short survey of models and of developed by me
methods using to optimize general networks with nonlinear non-convex
restrictions and objectives possessing mixed continuous-discrete
optimization variables. There are discussed the problem formulations and
solution methods for simulation, large-scale nonlinear discrete-continuous
network optimization, sensitivity and stability analysis and reliability
study for flow with nonlinear potential in general networks. The developed
methods belong to mixed integer nonnlinear programming (MINLP) in networks.
They use graph theory, integer and nonnlinear programming as their bases
corresponding to the nature of problems.

The mentioned above problems find an industrial application. The developed
methods are realized in software called ACCORD which is in operation in
different countries.

\section{Mathematical model and problem formulation}

\subsection{Passive and active elements}

A passive element is a two-pole that does not need an additional energy
supply to work. As simple examples, the resistors and pipe could be mention.
A dynamic model of a pipe could be written as a system of algebraic and
partial differential equations (PDE) for state variable vector $\left( \rho
,T,p,v,h\right) ^{T}:$%
\begin{equation}
\rho =\rho \left( p,h\right) :=\QDABOVE{1pt}{p}{zRT},\text{ where }z=z\left(
p,T\right)  \label{thermicStateEq}
\end{equation}
\begin{equation}
h=c_{p}T,\text{ where }c_{p}=c_{p}\left( p,T\right)
\label{calorifiicStateEq}
\end{equation}
\begin{equation}
\QDABOVE{1pt}{\partial }{\partial t}\left( 
\begin{array}{c}
\rho \\ 
\\ 
\rho v \\ 
\\ 
\rho \left( h+\QDABOVE{1pt}{v^{2}}{2}\right) -p%
\end{array}
\right) +\QDABOVE{1pt}{\partial }{\partial x}\left( 
\begin{array}{c}
\rho v \\ 
\\ 
\rho v^{2}+p \\ 
\\ 
\rho v\left( h+\QDABOVE{1pt}{v^{2}}{2}\right)%
\end{array}
\right) =\left( 
\begin{array}{c}
0 \\ 
\\ 
-\QDABOVE{1pt}{\lambda }{D}\rho \QDABOVE{1pt}{v\left| v\right| }{2}-\rho g%
\QDABOVE{1pt}{dH}{dx} \\ 
\\ 
\QDABOVE{1pt}{\Omega }{S}-v\rho g\QDABOVE{1pt}{dH}{dx}%
\end{array}
\right)  \label{dynHydro}
\end{equation}
where $c_{p}$ - specific heat by constant pressure $\left[ J/kg\ K\right] ;$ 
$\dot{M}$ - mass flowrate $\left[ kg/s\right] $; $g$ - gravitational
acceleration $\left[ m/s^{2}\right] $; $h$ - specific ehthalpy $\left[ J/kg%
\right] $; $F$ - objective function; $H$ - altitude $\left[ m\right] $; $p$
- pressure $\left[ N/m^{2}\right] $; $R$ - individual gas constant $\left[
J/kg\ K\right] $; $\rho $ - specific mass (density) $\left[ kg/m^{3}\right] $%
; $S$ - cross section $\left[ m^{2}\right] $; $T$ - temperature $\left[ K%
\right] $; $t$ - time coordinate $\left[ s\right] $; $v$ - flow velocity $%
\left[ m/s\right] $; $x$ - length coordinate $\left[ m\right] $; $z$ -
compressibility factor $\left[ -\right] $; $\Omega $ - heat flow per length $%
\left[ W/m\right] $. The equations express respectively thermic and
calorific state equations and conservation laws of mass, impulse, and
energy. The PDE is of hyperbolic type. For a practical simulation, a mixing
from initial and boundary conditions is used.

For a lot of considerations, a steady state model of a passive element is
precise enough. The model is algebraic and could be obtain after
simplification from (\ref{dynHydro}) as follows: 
\begin{equation}
q(x)=Const,\quad i\leq x\leq k  \label{qContin}
\end{equation}
\begin{equation}
f(p_{i},p_{k},q)=0,\quad \text{where}\quad p_{i}=p\mid _{x=i},\quad
p_{k}=p\mid _{x=k}  \label{tStead}
\end{equation}

As example of active elements, an amplifier, a pump, and a compressor could
be mentioned. Their common feature is that every of which have a plane
simply connected operational range \cite{Ot-98pat}. This range must be taken
as a restriction in the model. Unfortunately, the range forms a non-convex
domain mostly.

In some cases, in particular for hydraulic networks, active elements react
much quicker than passive elements. Then steady state models could be used
for active elements by network modelling. These steady state models are
algebraic.

Together with flow balance equations, the equations and inequalities
describing models of active and passive elements form a part of restrictions
by network modelling.

\subsection{Steady state continuous-discrete\ optimization problem in
general network}

Let be given a network $G=(V,E)$ with a set of nodes $V$ and a set of edges $%
E$. Let $q_{ik}$ is a flow in the edge $(i,k)$ and $Q_{i}$ is supply or
demand called intensity of node $i$. It means that 
\begin{equation}
\sum_{k\in V}q_{ik}+Q_{i}=0,\quad \text{ }i\in V,  \label{qQ}
\end{equation}
\begin{equation}
q_{ik}=-q_{ki},\quad (i,k)\in E.  \label{qarc}
\end{equation}
Suppose that there are given families of functional dependencies on flow and
potentials:

\begin{equation}
f_{d_{ik}}(p_{i},p_{k},q_{ik},c_{ik})=0,\quad (i,k)\in E;  \label{fdik}
\end{equation}%
\begin{equation}
d_{ik}\in \left\{ 1,...,N_{ik}\right\} ,\quad (i,k)\in E.  \label{dgr}
\end{equation}%
Here $c_{ik}$ is a vector of continuous parameters (coefficients), and $%
d_{ik}$ is a discrete parameter of the edge $(i,k).$ Suppose that there are
given limitations $Q_{i}^{-},$ $Q_{i}^{+},$ $p_{i}^{-},$ $p_{i}^{+},$ $%
c_{ik}^{-},$ $c_{ik}^{+}:$%
\begin{equation}
Q_{i}^{-}\leq Q_{i}\leq Q_{i}^{+},\quad i\in V;  \label{Qgr}
\end{equation}%
\begin{equation}
p_{i}^{-}\leq p_{i}\leq p_{i}^{+},\quad i\in V;  \label{pgr}
\end{equation}%
\begin{equation}
c_{ik}^{-}\leq c_{ik}\leq c_{ik}^{+},\quad (i,k)\in E;  \label{cgr}
\end{equation}%
and the other restrictions can be represented by inequalities with given $%
a_{ik}^{-},a_{ik}^{+}$ for the functions $%
a_{ik}(p_{i},p_{k},q_{ik},c_{ik},d_{ik})$ which have to be calculated:

\begin{equation}
a_{ik}^{-}\leq a_{ik}(p_{i},p_{k},q_{ik},c_{ik},d_{ik})\leq a_{ik}^{+},\quad
(i,k)\in E.  \label{agr}
\end{equation}%
The considering objective function depends both on flow $q_{ik}$ and on
potentials $p_{i},p_{k},$ on intensities $Q_{i},$ continuous $c_{ik}$ and
discrete parameters $d_{ik}.$ The problem consists in 
\begin{equation}
\text{minimize }F=\sum_{(i,k)\in
E}F_{ik}^{(1)}(p_{i},p_{k},q_{ik},c_{ik},d_{ik})+\sum_{i\in
V}F_{i}^{(2)}(p_{i},Q_{i})  \label{min}
\end{equation}%
\begin{equation}
\text{subject to (\ref{qQ} - \ref{agr}).}  \label{restrict}
\end{equation}%
The set of available values of discrete parameters $d_{ik}$ in (\ref{dgr})
means that a family of functions $f_{d_{ik}}$ can act on the edge $(i,k),$
and we have to select the best form of every equation (\ref{fdik}). Thus, we
have to select the best equation (\ref{fdik}) from the point of view of
functional (\ref{min}). We may interpret this as a selection of the most
profitable equipment which is installed or can be installed on the place of
the edge $(i,k)$.

The continuous parameters $c_{ik}$ in (\ref{fdik}) can be interpreted as
parameters that smoothly regulate the work of equipment $d_{ik}$ in bounds (%
\ref{cgr}). We may interpret (\ref{agr}) as restrictions on the power,
temperature, dissipation and other characteristics of the equipment that is
represented on the edge $(i,k)$ . The inequalities (\ref{Qgr}) and the
dependence of objective function on intensities $Q_{i}$ mean that the most
profitable values of supplies and demands have to be chosen.

The functions $F_{ik}^{(1)},f_{d_{ik}},a_{ik}$ staying in the objective
function and in the restrictions are non-linear and can be non-smooth and
non-convex but expensive to compute. Usually the number of nodes and edges
are hundreds with a tendency to be thousands. So the formulated problem is a
problem of large-scale nonlinear discrete-continuous programming in general
networks. This problem generalize the well-known minimal cost network flow
problem.

\section{Methods of the discrete-continuous nonlinear network optimization}

The formulated problem (\ref{qQ} - \ref{min}) is a problem for searching a
conditional extremum for a function with continuous-discrete parameters of
optimization. According to the mixing character of variables, the
combination of the continuous and discrete optimization methods is used for
the problem solving. The integer and nonlinear programming and graph theory
are basis of the offered method. Its main characteristic feature consists in
the obtaining of dominant solutions on fragments of the network.

\subsection{Simple examples}

Let discrete parameters be fixed. The next 2 examples and orders are dual to
each other.

\begin{example}
Let a node $v_{root}$ is chosen. There exists a linear order $\ell _{p}$ on $%
V$ such that if we enumerate the nodes with respect to this linear order
then: 1) the root has the number $n_{v_{root}}=1$; 2) for every node with
number $n>1$ there is just one adjacent node with number $n_{u}<n.$  By
other words, it is possible to construct an optimization process so that
among potentials $p_{i}$ only $p_{root}$ is considered as an independent
variable , and all other $p_{i}$ are dependent variables computing after $%
p_{root}$ and $q_{ik}.$
\end{example}

\begin{example}
There exists a linear order $\ell _{q}$ on nodes $V$ such that if we
enumerate the nodes with respect to this linear order then for every node
with number $n$ the numbers of all adjacent nodes except may be one are less
than $n:n_{u}<n$ . By other words, it is possible to construct an
optimization process so that among flows $q_{ik}$ only non-tree flows are
considered as independent variables, and all other flows $q_{ik}$ are
dependent variables.
\end{example}

In terminology of linear programming, these independent variables correspond
to non-basic variables, and .dependent variables are basic ones.

\subsection{The continuous dominant problems}

Let $D$ and $X$ are correspondingly the sets of discrete and continuous
independent variables. Denote a problem for searching a feasible solution (%
\ref{qQ} - \ref{agr}) by $P(D,X)$ and an optimization problem (\ref{qQ} - %
\ref{min}) by $FP(D,X)$.

Problem $P\symbol{94}$ dominates another one $P$ and is called as a dominant
for $P$, if the feasible set for $P\symbol{94}$ contains such set for $P$.

\begin{lemma}
\label{dominant}For any discrete-continuous problem $P(D,X)$ there exists a
continuous dominant $P\symbol{94}(X_{D},X)$, where a set of discrete
variables $D$ is replaced by a set of continuous variables $X_{D}.$
\end{lemma}

\proof
On the edges with a discrete variable $d_{ik},$ the equation $%
t_{ik}=p_{i}-p_{k}$ replaces the family of equations (\ref{fdik} - \ref{dgr}%
). A component $t_{ik}$ of tension vector replaces the discrete variable $%
d_{ik}$ and the continuous ones $c_{ik}$ there. 
\endproof%

As consequence, there is the continuous dominant $FP\symbol{94}(X_{D},X)$
for an initial optimization problem $FP(D,X).$

Assign a node where potential is either given or varied as a root. Let a set
of discrete variables $D$ is participated onto two classes $D^{\prime }$ and 
$D^{C}$. Collections $d^{\prime }\in D^{\prime }$ are called as fragments of 
$d\in D.$ Consider an initial discrete-continuous problem (\ref{qQ} - \ref%
{agr}) by fixed $d^{\prime }\in D^{\prime }$; denote it $P_{d^{\prime
}}(D^{C},X)$. We can construct for a problem $P_{d^{\prime }}(D^{C},X)$ its
continuous dominant $P_{d^{\prime }}\symbol{94}(X_{D^{C}},X)$. The problem $%
P_{d^{\prime }}\symbol{94}$ defines a feasibility of fragment $d^{\prime }$.
Fragment $d^{\prime }$ is infeasible, if $P_{d^{\prime }}\symbol{94}$ has no
solution. The explicit definition of fragment $d^{\prime }$ is based on the
follows

\begin{lemma}
\label{M-fragment}Let $G=(V,E)$ is a connected graph, $v_{0}\in V$ is a
node, and a finite set $M_{0}\subset E$ consists of $\left| M_{0}\right| $
edges. Consider $M$ as a set of subsets of $M_{0}$ and write down $d\in M$
as 
\begin{equation}
\begin{array}{lll}
d=(d_{1},...,d_{\left| M_{0}\right| })_{\rho }\quad ,\quad  &  & d_{i}=0,1.
\end{array}
\label{M-present}
\end{equation}
for a linear order $\rho $ on $M_{0}.$ There exists such a linear order $\pi 
$ on $M_{0}$ that for every $m$-length fragment $%
d^{m}=(d_{1},...,d_{m},0,...,0)_{\pi }$ there is a connected subgraph $%
G_{i}=(V_{i},E_{i}),i=i(m),$ that is a minimal connected subgraph containing
the node $v_{0},$ all the edges from $d^{m},$ and all their endnodes.
\end{lemma}

\proof
A proof is constructive. It consists of the following steps:

1) Construction a linear order $\ell _{v}$ on a rooted spanning tree;

2) Construction a linear order $\pi $ on the subset $M_{0}$ of edges;

3) Construction a minimal subgraph for a set $d^{m}$ of the first $m$ edges
from $M_{0}$ . 
\endproof%

Let the discrete variables are acting only on edges. Take the set of edges
possessing discrete variables in the role of $M_{0}$. Let an edge $e\in
M_{0} $ has a number $j$ in a linear order $\pi ;$ denote it as $e_{j}.$
With the notation (\ref{M-present}) for $d\subseteq M_{0}$ and $%
d^{m}=(d_{1},...,d_{m},0,...,0)$ for $m$ - length fragment, let us allow to
use a non-null value of a discrete variable acting on the edge $e_{j}\in d$
as $d_{j},$ and to use $d_{j}=0$ for $e_{j}\notin d$ . In the last case we
shall say that discrete variable is not given there. So the linear order $%
\pi $ on $M_{0}$ induces a linear order on a set $D$ of discrete variables,
and we can tell about $m$-length fragment of discrete variables. The lemma %
\ref{M-fragment} can be reformulated as follows.

\begin{theorem}
\label{fragment}Let $G=(V,E,v_{root})$ is connected network with one node
marked as $v_{root}$. Let some edges of the network possess discrete
variables, and the last ones are only there. There exists a linear order on
the set of discrete variables $D$ that:

1) every $d\in D$ can be written as 
\begin{equation}
\begin{array}{ll}
d=(d_{1},...,d_{M}), & d_{i}=0,1,...,N_{i},
\end{array}
\label{d-full}
\end{equation}
where $d_{i}=0$ means that value of corresponding variable is not given;

2) for every $m$-length fragment $d^{m}=(d_{1},...,d_{m},0,...,0)$ there
exists a connected sub-network $G_{i}=(V_{i},E_{i}),i=i(m),$ that is a
minimal sub-network containing the root, all the edges with non-null
discrete variables from $d^{m},$ and all endnodes of these edges.
\end{theorem}

The solving method for the continuous-discrete problem $FP(D,X)$ can be
presented as a branching multi-level computational process. On the main
branch of this process we solve the continuous problem $FP\symbol{94}(D,X)$
or $P(D,X)$.

\subsection{A method for solving the continuous problems in general networks}

In continuous problems $FP(X),FP\symbol{94}(X_{D},X),$ and $P_{d^{\prime }}%
\symbol{94}$ , the components of continuous variable vector $%
x=(x_{1},...,x_{n})$ are divided into 3 groups: those from the initial
problem $FP(D,X)$; the flows in chords; and the continuous variables which
change the discrete ones in the dominant continuous search. The flows in
chords are of a special meaning. They can be computed as a solution of the
equation system (\ref{qQ}),(\ref{fdik}),(\ref{dgr}) by the fixed both
discrete and continuous variables of 2 other groups. Then the dimension of
continuous variables and of equation systems (\ref{qQ}),(\ref{fdik}) can be
reduced thanks to the direct solving of such systems on a tree.

To solve the continuous search problem, the classical methods of constrained
and unconstrained optimization from nonlinear programming are implemented.
The preference is given to the methods without derivatives because of
non-convex and non-smooth functions in restrictions and objectives.

Components of objective, penalty and Lagrangian functions are computed as
soon as a network state component $s_i$ is been known, $%
s=(q_{ik},Q_i,p_i,c_{ik},d_{ik}).$

By computation of state $s$, we have to distribute flow and potential. To
effectively solve the equation systems (\ref{qQ}) and (\ref{fdik}) on tree,
respectively the width- and depth-first searches in graph are used.

\subsection{Modification of the branch-and-bound method}

In our modification of the branch-and-bound method, only $m$-length
fragments are used. If fragment is infeasible then all its extension $%
(d_{1},...,d_{m},j_{m+1},...,j_{M})$ are excluded from the further
examination. The examination is gone in a lexicographic order of vectors $%
d=(d_{1},...,d_{M})$. Starting from $d^{(0)}=(0,...,0),$ a transition from $%
d^{(n)}$ to $d^{(n+1)}$ is performed by the next rule: if $d^{(n)}=d^{m}$ is
fragment with length $m<M$ and $d^{(n)}$ is feasible then 
\begin{equation}
d^{(n+1)}=(d_{1},...,d_{m},1,0,...,0)  \label{lengthen}
\end{equation}%
else $d^{(n+1)}=d^{m}\bigoplus e^{m}.$ Here $e^{m}=(0,...,1,...0)$ and 
\begin{equation}
d^{m}\bigoplus e^{m}=\left\{ 
\begin{array}{ll}
(d_{1},...,d_{m-1},d_{m}+1,0,...,0), & \text{if }d_{m}<N_{m}; \\ 
(d_{1},...,d_{\mu -1}+1,0,...,0), & 
\begin{array}{l}
\text{if }d_{m}=N_{m}\text{ and } \\ 
\exists \mu :1<\mu \leq m:%
\end{array}
\\ 
& 
\begin{array}{l}
d_{\mu -1}<N_{\mu -1}\text{ \& }d_{k}=N_{k},\text{ } \\ 
k=\mu \ldots m;%
\end{array}
\\ 
(0,...,0), & \text{if }d_{k}=N_{k},\text{ }\forall k=1,\ldots ,M.%
\end{array}%
\right.  \label{bab}
\end{equation}

\subsection{Another procedure for searching an optimal integer-feasible
solution}

The idea of another procedure for $FP(D,X)$ is connected with a possibility
to check, if there is a realization of equation $t_{ik}=p_{i}-p_{k}$ in the
equation family (\ref{fdik}),(\ref{dgr}). This sub-problem is to construct
mapping 
\begin{equation}
\left\{ \left( t_{ik},p_{i},p_{k},q_{ik}\right) \mid
t_{ik}=p_{i}-p_{k}\right\} \longrightarrow \left\{ \left(
d_{ik},c_{ik},p_{i},p_{k},,q_{ik}\right) \mid (\ref{fdik}),(\ref{dgr}),(\ref%
{cgr})\right\} .  \label{td}
\end{equation}%
The second procedure for solving $FP(D,X)$ is proposed as follows. The
dominant problem $FP\symbol{94}(X_{D},X)$ is to be solved. On the every step
of continuous optimization $FP\symbol{94}(X_{D},X)$, there is checked the
existence of feasible discrete variable $d_{ik}$ and perhaps continuous ones 
$c_{ik}$ for every continuous variable $t_{ik}$ changing $d_{ik}$ in the
dominant problem. If there is no such feasible discrete variable $d_{ik}$
then the corresponding continuous variable $t_{ik}$ would be penalized. In
other words, $t_{ik}$ is penalized if there is no mapping (\ref{td}). Then
an optimum of problem $FP\symbol{94}(X_{D},X)$ brings an optimum for the
initial problem $FP(D,X).$

The difference between this method and the branch-and-bound approach is
clear for the cases of pure discrete problems on a tree. Such class of
problems means that graph $G=(V,E)$ is a tree and there is no continuous
variables at all: $X=\O .$ Then the formulated version of branch-and-bound
method can be very efficient. It brings a global optimum as a solution.

For the same problem $FP(D)$ on a tree, the second method based on the
checking of integer-feasibility by the solving of dominant problem $FP%
\symbol{94}(X_{D})$ can fall in a local optimum. To avoid it, a restart can
be used. By restart, it is worth to move the root into another node, if the
restrictions (\ref{pgr}) allow it.

However, for the network with complicated topology the branch-and-bound
method could be expensive if not prohibited. This is valid especially for
the networks with a lot of cycles. Then the second procedure with penalizing
non-realizable $t_{ik}$ is available there. In practice, the method avoids
the local optima.

\section{Risk, sensitivity and stability analysis for reliability study}

There are 2 extreme approaches by analysis either a network state is stable.
The first one is used by operating when a control vector is chosen (fixed)
.and it is necessary to check either this control vector still lead to a
feasible network state if initial or boundary conditions are changing. The
second one is used by design when it is necessary \ to check either there is
a value of a control vector for every possible operating point. The first
approach we shall call the strong stability while the second one we shall
call the weak stability.

Let be given a value $\pi _{0}$ of analyzing parameter $\pi ,$ a fixed
control $u_{0},$ and $\delta _{0}>0.$ The scalar parameter $\pi $ is called 
\textbf{feasible strong stable }in the point $\pi _{0}$ by the fixed control 
$u_{0}$ if there is an interval $\left( \pi _{0}-\delta ,\pi _{0}+\delta
\right) ,$ $\delta >0,$ such that the network state $x(u_{0},\pi )$ is
feasible for every $\pi \in \left( \pi _{0}-\delta ,\pi _{0}+\delta \right) $%
. If $\delta \geq \delta _{0}$ then the network state $x_{0}=x(u_{0},\pi
_{0})$ is called feasible strong stable with respect to the scalar parameter 
$\pi $ by the fixed control $u_{0}.$

To prove the stability of parameter $\pi $ in the neighborhood of $\pi _{0}$
, we find such interval $\left[ \pi _{1},\pi _{2}\right] $ containing $\pi
_{0}$ that the network state $x(u_{0},\pi )$ is feasible for every $\pi \in %
\left[ \pi _{1},\pi _{2}\right] $ , and for every small $\varepsilon >0$
there are $\pi _{1}\in \left[ \pi _{1}-\varepsilon ,\pi _{1}\right] ,$ $\pi
_{2}\in \left[ \pi _{2},\pi _{2}+\varepsilon \right] $ that for $\pi =\pi
_{1},\pi _{2}$ the network state $x(u_{0},\pi )$ is not feasible.

Let be given a value $\pi _{0}$ of analyzing parameter $\pi ,$ a
neighborhood $U_{0}$ of a control $u_{0}$ and $\delta _{0}>0,$ $\eta >0.$
Let $F\left( u\right) =F\left( u,x\left( u,\pi \right) \right) $ be an
objective function. The scalar parameter $\pi $ is called \textbf{feasible
weak stable }in the point $\pi _{0}$ if there is an interval $\left( \pi
_{0}-\delta ,\pi _{0}+\delta \right) ,$ $\delta >0,$ such that for every $%
\pi \in \left( \pi _{0}-\delta ,\pi _{0}+\delta \right) $ there is a control 
$u\in U_{0}$ which provides that the network state $x(u,\pi )$ is feasible
and that 
\begin{equation}
\left| F\left( u_{0},x(u_{0},\pi \right) -\min_{u\in U_{0}}F_{1}\left(
u,x\left( u,\pi \right) \right) \right| <\eta  \label{weakStabl}
\end{equation}
where $F_{1}\left( u,x\right) =F\left( u,x\right) +\psi \left( u\right) $
and $\psi \left( u\right) $ is a cost of control switch from $u_{0}$ to $u.$

For a scalar parameter $\pi ,$ it is possible to compute the above interval $%
\left[ \pi _{1},\pi _{2}\right] $ showing a complete feasibility interval
with respect to parameter $\pi .$ For a set of analyzing parameters $\pi
=\left( \pi _{1},...,\pi _{n}\right) ,$ a statistic test method could be
used. If there are $\delta _{i}>0$ that for the most points $\pi \in \Pi
\left( \pi _{i0}-\delta _{i},\pi _{i0}+\delta _{i}\right) $ (saying, for 95
\%) there is a control $u\in U_{0}$ providing feasibility of the network
state $x(u,\pi )$ and condition (\ref{weakStabl}) then network state is
called weak stable with respect to parameter set $\pi .$ \ The definition of
strong stability with respect to parameter set $\pi $ is evident.

In complement to statistic approach, one simple result for 'definite'
approach has to be mentioned. The capacity $q_{ik}^{+}$ of an edge $(i,k)$
can be found from equation (\ref{fdik}) varying $p_{i},p_{k},d_{ik}$ : 
\begin{equation}
q_{ik}^{+}:=\max_{d_{ik},\ p_{i},p_{k}}\left\{ 
\begin{array}{c}
q_{ik}\mid f_{d_{ik}}(p_{i},p_{k},q_{ik})=0,\quad \\ 
d_{ik}\in \left\{ 1,...,N_{ik}\right\} ,\quad \\ 
p_{n}\in \left[ p_{n}^{-},p_{n}^{+}\right] ,\quad n=i,k%
\end{array}%
\right\}  \label{defCap}
\end{equation}%
In some cases we can precise potential intervals $\left[ p_{i}^{-},p_{i}^{+}%
\right] $ in (\ref{pgr}) and capacities therefore.

\begin{lemma}
Let be given intervals $\left[ p_{i}^{-},p_{i}^{+}\right] $ (\ref{pgr}),
functions $f_{d_{ik}}$ in models (\ref{fdik}) are monotonous relative both
potentials $p_{i},p_{k}$ and $q_{ik}$ separately, and every model $f_{d_{ik}}
$ in (\ref{fdik}) is chosen: $d_{ik}=Const,\quad (i,k)\in E,$ i.e. $%
f_{d_{ik}}=f_{_{ik}}$.  For a known flow $q_{ik}$ , it can be found such
intervals 
\begin{equation}
\left[ p_{i}^{\min },p_{i}^{\max }\right] \subseteq \left[
p_{i}^{-},p_{i}^{+}\right]   \label{preciseP}
\end{equation}
that the inequality system (\ref{fdik}),(\ref{pgr}) has a solution if and
only if the system (\ref{fdik}),(\ref{preciseP}) has a solution.
\end{lemma}

In other words, for a node $j\in V$ every $p_{j}^{0}\notin \left[
p_{j}^{\min },p_{j}^{\max }\right] $ leads to a solution of equation system $%
\left\{ (\ref{fdik}),p_{j}=p_{j}^{0}\right\} $ which violates the inequality
system (\ref{fdik}), (\ref{pgr}).

\section{Conclusions}

The solution methods and their implementation are described for a class of
problems of large-scale nonlinear discrete-continuous network optimization
in general networks. Objective function depends both on flow and potential.
Minimal cost and maximal flow problems are generalized. The problem consists
both in the optimal choice of dependencies on flow and potential from the
families given on some edges and in optimization of values of intensities,
flow, and potential.

The formulations, methods and algorithms are developed for the following
problems: 1) the (pure) continuous nonlinear optimization in general network
with a given functional dependence between flow and potentials; 2) the
(pure) discrete nonlinear optimization problem of selection the best
functional dependence between flow and potentials from given families; 3)
the (mixed) continuous-discrete nonlinear optimization in general network
with given families of dependences between flow and potentials. The
developed optimization method represents a branching multi-level
computational process. It is based on nonlinear and integer programming and
on graph theory. Its main characteristic feature consists in obtaining of
dominant solutions on network fragments. During optimization, the offered
methods adapt themselves to the structure of every network and to the
features of both objective function and variables. It allows to use such
objectives as minimization of cost, of expenditure, and of set-point
deviation, maximization of flow, of profit, and so on.

On the base of the developed optimization methods, some approaches for
sensitivity analysis, stability investigation and reliability study are
proposed.

\end{document}